\documentclass{amsart}

\usepackage{amssymb,latexsym,graphicx}

\newtheorem{theorem}{Theorem}[section]
\newtheorem{lemma}[theorem]{Lemma}

\theoremstyle{definition}
\newtheorem{definition}[theorem]{Definition}

\theoremstyle{remark}
\newtheorem{remark}[theorem]{Remark}

\numberwithin{equation}{section}

\begin{document}

\title[Locating the zeros of partial sums of $e\sp z$ with Riemann-Hilbert]
{Locating the zeros of partial sums of $e\sp z$ with Riemann-Hilbert methods}

\author{T. Kriecherbauer}
\address{Fakult\"at f\"ur Mathematik, Universit\"at Bochum, Universit\"atsstr. 150,
D-44801 Bochum, Germany}
\email{thomas.kriecherbauer@rub.de}
\thanks{The first and second author are supported by the
Belgian Interuniversity Attraction Pole P06/02.}
\thanks{The first and third author were supported in part by the SFB/TR 12 of the
Deutsche Forschungsgemeinschaft.}

\author{A. B. J. Kuijlaars}
\address{Department of Mathematics, Katholieke Universiteit Leuven, Celestijnenlaan 200 B, 3001 Leuven, Belgium}
\email{arno.kuijlaars@wis.kuleuven.be}
\thanks{The second author is supported by FWO-Flanders project G.0455.04,
by K.U. Leuven research grant OT/04/21, by the
 European Science Foundation Program MISGAM,
and by a grant from the Ministry of Education and
Science of Spain, project code MTM2005-08648-C02-01.}

\author{K. D. T-R McLaughlin}
\address{Department of Mathematics, University of Arizona, Tucson, AZ 85721, USA}
\email{mcl@math.arizona.edu}
\thanks{The third author was partially supported by the NSF under
grants DMS-0200749 and DMS-0451495.}

\author{P. D. Miller}
\address{Department of Mathematics, University of Michigan, 530 Church St.,
Ann Arbor, MI 48109, USA}
\email{millerpd@umich.edu}
\thanks{The fourth author was partially supported by the NSF under grants
DMS-0103909 and DMS-0354373, and by a grant from the Alfred P. Sloan foundation.}


\dedicatory{Dedicated to Percy Deift with gratitude and admiration.}


\begin{abstract}
In this paper we derive uniform asymptotic expansions for the
partial sums of the exponential series. We indicate how this
information will be used in a later publication to obtain
full and explicitly computable asymptotic expansions with error bounds
for all zeros of the Taylor polynomials $p_{n-1}(z) =
\sum_{k=0}^{n-1} z^k/\,k!\,$. Our proof is based on a
representation of $p_{n-1}(nz)$ in terms of an integral of
the form $\int_{\gamma} \frac{e^{n\phi(s)}}{s-z}ds$. We
demonstrate how to derive uniform expansions for such
integrals using a Riemann-Hilbert approach. A comparison with
classical steepest descent analysis shows the advantages of the
Riemann-Hilbert analysis in particular for points $z$ that are close to the
critical points of $\phi$.
\end{abstract}

\maketitle

\section{Introduction}
\label{S1}

During the past fifteen years and largely due to the ground breaking
work \cite{DZ1}, \cite{DZ2} of Deift and Zhou, Riemann-Hilbert problems
have become a powerful tool in asymptotic analysis with
applications in many fields such as
inverse scattering theory, integrable PDE's, orthogonal polynomials,
statistical mechanics and random matrix theory. With this paper we would like
to add another item to this list of applications,
namely the study of the asymptotic behavior of zeros of Taylor polynomials
of entire functions. Although our method works in principle for a large class
of entire functions we will restrict our attention to the classic case of
Taylor polynomials of $e^z$ in order to keep the
presentation as non technical as possible. It is also one of the intentions of this paper
to highlight the advantageous aspects of the Riemann-Hilbert approach in a
simple situation. The first and maybe most surprising feature is how many
objects of interest can be characterized as the unique solution of some
Riemann-Hilbert problem and with this paper we present one more illustration of this fact.
Secondly, there is a set of techniques to transform Riemann-Hilbert problems
to some form where an asymptotic expansion of the quantity of interest can be extracted
in a systematic fashion. In our case, we will only need one of these techniques and that is
the construction of a local parametrix. This construction exemplifies in a very simple situation
just one of the wonderful ideas of Deift and Zhou \cite{DZ1}, \cite{DZ2}, namely how to effectively
localize a Riemann-Hilbert problem which by definition is a global problem in the complex plane.
Localization is an indispensable step in the asymptotic analysis of any Riemann-Hilbert problem where the
asymptotic behavior of the solution is different for different regions of the complex plane.
Moreover, and that is the third aspect we want to highlight, the localization procedure
of Deift and Zhou always provides for matching asymptotics in those regions where the
asymptotic behavior changes and one obtains error bounds uniform in all of $\mathbb C$.

We denote by $p_{n}(z) := 1 + z + \cdots + \frac{z^{n}}{n!}$ the partial sums of the
exponential series. The problem to describe the asymptotic distribution of the zeros
of $p_n$ was posed and solved in the classical paper of Szeg\H{o} \cite{S}. He proved
that the zeros of $p_n$, divided by $n$, converge in the limit $n \to \infty$ to
some curve $D_{\infty}$, now called
Szeg\H{o} curve, which consists of all
complex numbers $|z| \leq 1$ that satisfy the equation $|z e^{1-z}| = 1$.
\begin{figure}
\includegraphics{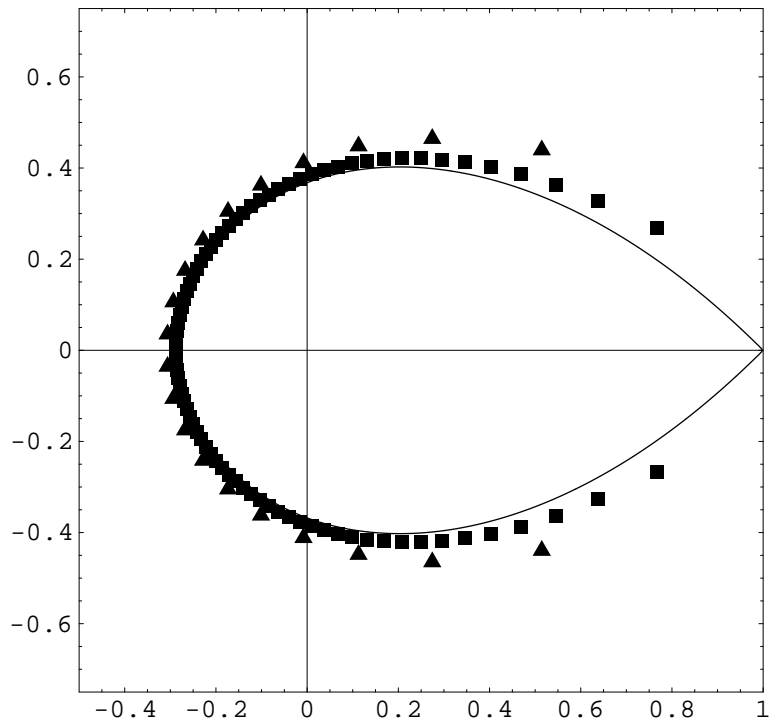}
\caption{Szeg\H{o} curve (solid line) and rescaled zeros of $p_{20}$ (triangles)
and $p_{80}$ (boxes)}
\label{f1.1}
\end{figure}
The curve $D_{\infty}$ together with the rescaled zeros of $p_n$ for $n=20$ and
$n=80$ are displayed in Figure \ref{f1.1}.
Moreover, Szeg\H{o} also determined the limiting distribution of the rescaled zeros
on $D_{\infty}$.
A first result presenting error bounds on the distance between the Szeg\H{o} curve and the
zeros of $p_{n}(nz)$ has been established by Buckholtz \cite{B} who showed that
they are located in the exterior of $D_{\infty}$ at a distance of at most $2e/\sqrt{n}$.
For each zero of $p_{n}(nz)$ its distance from $D_{\infty}$ is measured by the minimal distance
between the zero and all the points of the Szeg\H{o} curve.
Subsequently more detailed asymptotics of $p_{n}(nz)$ have been derived. It turns out
that $z_0=1$ is a critical point where the asymptotic behavior changes. It was shown
by Newman and Rivlin \cite{NR} that in neighborhoods of $z_0$ of size $O(1/\sqrt{n})$ the
rescaled polynomials $p_{n}(n+w\sqrt{n})$ can asymptotically be expressed in terms of the
complementary error function erfc.
Carpenter, Varga, and Waldvogel \cite{CVW} (see also \cite{VC1} and see \cite{V} for an interesting
discussion of zeros of the partial sums of $e^z$) then provided an asymptotic
expansion of $p_{n}(nz)$ in compact subsets of $\mathbb C \setminus \{z_0\}$. These
results were used in \cite{CVW} to obtain lower and upper bounds on the distance of
the zeros from the Szeg\H{o} curve. Note that up to the recent paper by Bleher and Mallison
\cite{BM} no uniform asymptotics for $p_{n}(nz)$ and its zeros were available in a
fixed size neighborhood of the critical point $z_0$.
In a different direction and
using methods from logarithmic
potential theory Andrievskii, Carpenter and Varga \cite{ACV} recently extended the
results of Szeg\H{o} \cite{S} on the angular distribution of the zeros by proving
uniform error bounds for regions including $z_0$.
We refer the reader to the recent review \cite{Z}
for a description of more results on the zeros of $p_n$ . Related results on the
zeros of the partial sums of $\cos$, $\sin$ and of more general sums of exponential
functions can be found in \cite{BM} and in references therein. Finally, we mention
the review \cite{O} on the behavior of the zeros of more general sections and tails of
power series.

The main novel result proved in the present paper is Theorem \ref{th3.1}
which provides an explicitly computable asymptotic expansion near $z_0$
for a quantity $F_n(z)$ closely related to $p_{n-1}(nz)$ (see (\ref{eq1.20}),
(\ref{eq1.40}) below). Together with Remark \ref{r3.2} we produce
an asymptotic expansion of $F_n$ with uniform
error bounds for $z$ contained
in some region $V$ where all the roots of  of $p_{n-1}(nz)$ are located.
These results will be used in a subsequent
publication to derive explicitly computable asymptotic expansions for all the
zeros of $p_{n-1}$ in terms of the zeros of the complementary error function
for zeros close to the critical point,
otherwise in terms of the solutions of $(z e^{1-z})^n = 1$ which lie on the Szeg\H{o}
curve. We will state the corresponding results without proof in
Theorems \ref{th4.2} and \ref{th4.1} below.

Let $\gamma$ be any smooth Jordan curve
encircling the origin counterclockwise.  For $z \in \mathbb C \setminus \gamma$ define
\begin{equation}
\label{eq1.20}
F_n(z) := \frac{1}{2 \pi i} \int_{\gamma}\frac{s^{-n}e^{n(s-1)}}{s-z}ds =
\frac{1}{2 \pi i} \int_{\gamma}\frac{e^{n \phi (s)}}{s-z}ds\, , \,
\mbox{ with } \phi(s) = s-1-\ln s \, ,
\end{equation}
where $\ln$ denotes the standard branch of the logarithm slit along the negative real axis.
A straightforward application of the calculus of residues then yields
\begin{equation}
\label{eq1.40}
(ez)^{-n} p_{n-1}(nz) =
\left\{ \begin{array}{ll}
-F_n(z)\, , &\mbox{ for $z$ in the exterior of $\gamma$ }\\
e^{n \phi (z)} - F_n(z)\, , &\mbox{ for $z \neq 0$ in the interior of $\gamma$ }
\end{array}
\right.
\end{equation}
To the best of our knowledge the Cauchy-type representation of $p_{n-1}$ provided by
(\ref{eq1.20}) and (\ref{eq1.40}) which is crucial for our Riemann-Hilbert approach
has so far not been used in the asymptotic analysis of the zeros of $p_{n-1}$.

Given the form of $F_n$ it is natural to use the method of steepest descent for the
asymptotic analysis of the integral. This method requires that the path $\gamma$ of
integration passes through
the critical point $z_0=1$ of the function $\phi$. This indicates already that the asymptotic
analysis will be most difficult for $z$ close to $1$ due to the term $(s-z)^{-1}$ in the
integrand. We discuss this approach in detail in Section \ref{S2} below and we will see
that by standard techniques of steepest descent analysis one may derive an asymptotic expansion
for $F_n(z)$ with uniform error bounds for $|z-1| > n^{- \alpha}$ and $0 < \alpha < 1/2$.

Section \ref{S3} contains the central result of this paper where we derive an asymptotic
expansion for $F_n$ in some fixed neighborhood $U_{\epsilon}(1)$ of the critical point. Here we
use that $F_n$ can easily be characterized as the unique solution of a Riemann-Hilbert problem
since by definition $F_n$ is just the Cauchy transform of $e^{n \phi}$ (see ($RHP$)$_1$ described
at the beginning of Section \ref{S3}). Using
the standard Riemann-Hilbert technique of constructing a local parametrix near the critical point
the derivation of the asymptotic expansion is no more difficult than the computations in
the non critical situation of Section \ref{S2}.

In the final Section \ref{S4} we briefly explain how to use the asymptotic results to
obtain asymptotic information on the location of the zeros. This, however, is a somewhat
technical affair and we will present a complete account of all results and their proofs in a
later publication. In Theorems \ref{th4.2}, \ref{th4.1} we state our results for the zeros of $p_{n-1}$
in the upper half plane.

\section{Classical steepest descent analysis}
\label{S2}

In order to remind ourselves of the method of steepest descent (see also \cite{Miller}
for an elementary exposition) we first
determine the large $n$ asymptotics of a quantity that is related to $F_n$
defined in (\ref{eq1.20}) but is somewhat simpler to analyze. Let
\begin{equation}
\label{eq2.10}
G_n : =  \frac{1}{2 \pi i} \int_{\gamma} e^{n \phi (s)} ds=
\frac{1}{2 \pi i} \int_{\gamma} s^{-n} e^{n(s-1)} ds \, ,
\end{equation}
where $\gamma$ denotes some smooth Jordan curve that is oriented
counterclockwise containing the origin in its interior.

\begin{figure}
\includegraphics{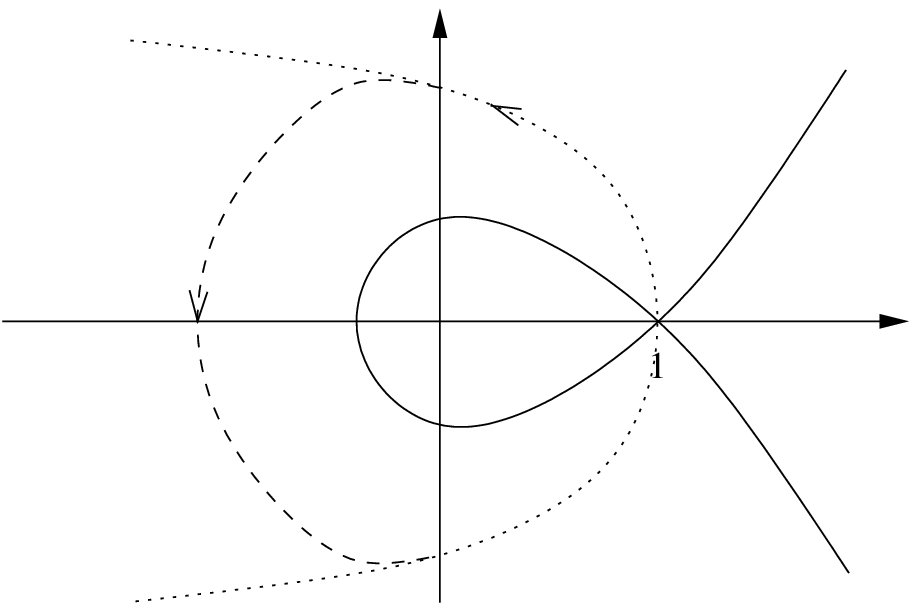}
\caption{}
\label{f2.1}
\label{p10}
\end{figure}

The method of steepest descent offers a recipe how to deform the
contour of integration in such a way that the asymptotic behavior
of $G_n$ can be determined most conveniently. More precisely, we are
to consider contours that pass through critical points of $\phi$
along the path of steepest descent with respect to the real part of
$\phi$. Such a choice of the contour ensures that the modulus of
the integrand obtains its local maximums only at the critical points
of $\phi$ and -- up to exponentially small error terms -- the integral
is determined by the contributions of those parts of the contour which
lie in small neighborhoods of the critical points. In the situation
at hand the function $\phi$ has only one critical point at $z_0 = 1$
and we can easily understand the behavior of the real part of $\phi$.
It has a saddle point at $z_0 = 1$ and the solid curve in Figure
\ref{f2.1} displays those numbers $z$ satisfying $\Re (\phi (z)) =
\Re (\phi (z_0))=0$. Note, that the Szeg\H{o} curve $D_{\infty}$ which
was displayed in Figure \ref{f1.1} above coincides with the closed loop
part of the solid curve.
Moreover, the dotted line in this figure
shows the path of steepest descent
away from the saddle. This path can easily be determined to consist of those
points in $\mathbb C \setminus \{0\}$ satisfying $\arg z=\Im z$ by using
the property that the
imaginary part of $\phi$ remains constant on it, i.e.
$\Im (\phi (z)) = \Im (\phi (z_0))= 0$. Choosing the smooth closed contour
$\tilde{\gamma}$ to coincide with the path of steepest descent in the
right half plane and with the dashed line in the left half plane
(in the left half plane any smooth curve that does not intersect the
solid line will do) we obtain
\begin{equation}
\label{eq2.20}
G_n = \frac{1}{2 \pi i} \int_{\tilde{\gamma}} e^{n \phi (s)} ds =
\frac{1}{2 \pi i} \int_{\tilde{\gamma}\cap U} e^{n \phi (s)} ds
+ O\left(e^{-nc} \right) \, ,
\end{equation}
where $U$ is any fixed neighborhood of $z_0=1$. It is
only the number $c > 0$ in the error term that will depend on the
choice of $U$.

In order to further analyze the integral along $\tilde{\gamma}\cap U$
we change variables near the critical point of $\phi$.
Since $\phi^{\prime\prime}(z_0) \neq 0$ we can find an open neighborhood
$U_0$ of $z_0$ and a biholomorphic map $\lambda: U_{\delta_0}(0)
\to U_0$ for some $\delta_0 > 0$ such that
\begin{equation}
\label{eq2.30}
\phi (\lambda (\xi)) = \xi^2 \quad \mbox{ for } \xi \in U_{\delta_0}(0)
\end{equation}
Observe that $\lambda$ maps the imaginary axis onto the path
of steepest descent, more precisely $\lambda (i \mathbb R \cap U_{\delta_0}(0))
= \tilde{\gamma} \cap U_0$. This can be seen from the
characterization of the path of steepest descent by the imaginary
part of $\phi$, i.e.
$$
\tilde{\gamma} \cap U_0 = \{z \in U_0 \colon
\Im \phi(z) = \Im \phi(1) =0 \mbox{ and } \Re \phi(z) \leq \Re \phi(1)= 0 \} \, .
$$
Choosing $U:=U_0$
in (\ref{eq2.20}) we obtain
\begin{equation*}
G_n = \frac{1}{2 \pi i} \int_{-i \delta_0}^{i \delta_0} e^{n \xi^2}
\lambda^{\prime}(\xi) d\xi + O\left(e^{-nc} \right)
= \frac{1}{2 \pi \sqrt{n}} \int_{- \delta_0 \sqrt{n}}^{\delta_0\sqrt{n}}
e^{-t^2} \lambda^{\prime}(it/\sqrt{n}) dt + O\left(e^{-nc} \right)
\end{equation*}

Expanding $\lambda$ at $\xi=0$ leads to an asymptotic expansion for
$G_n$. For example, we learn from
$$
\lambda(\xi) = 1 + \sqrt{2}\xi +\frac{2}{3}\xi^2 + \frac{\sqrt{2}}{18}
\xi^3 + O(\xi^4)
$$
that
$$
G_n =  \frac{1}{\sqrt{2 \pi n}} \left( 1 - \frac{1}{12}n^{-1} +
O(n^{-2}) \right).
$$
Observe that the calculus of residues implies $G_n= e^{-n} n^{n-1}/\Gamma (n)$
directly from (\ref{eq2.10}) and we have thus found some version of
Stirling's formula.

It is clear that we may extend our reasoning to integrals of the form
\begin{equation}
\label{eq2.40}
\frac{1}{2 \pi i} \int_{\gamma} h(s) e^{n \phi (s)} ds
\end{equation}
for analytic functions $h$. Up to exponentially small error terms
(\ref{eq2.40}) is given by
$$
\frac{1}{2 \pi \sqrt{n}} \int_{- \delta_0 \sqrt{n}}^{\delta_0\sqrt{n}}
h(\lambda(it/\sqrt{n}))e^{-t^2} \lambda^{\prime}(it/\sqrt{n}) dt.
$$
This expression can easily be expanded in powers of $1/n$ using the
expansions of $\lambda$ around $0$ and of $h$ around $1$. For
example, in the case $h(1) \neq 0$ we immediately obtain
\begin{equation}
\label{eq2.50}
\frac{1}{2 \pi i} \int_{\gamma} h(s) e^{n \phi (s)} ds =
\frac{h(1)}{\sqrt{2 \pi n}} \left( 1 + O (n^{-1}) \right)
\end{equation}
for any Jordan curve winding counterclockwise around
the origin.

We are now in a position to apply our reasoning to
$F_n$ by replacing the analytic function
$h$ in (\ref{eq2.50}) by the meromorphic functions
$h_z(s) = (s-z)^{-1}$. Indeed, we have
$$
F_n(z) = \frac{1}{2 \pi i} \int_{\gamma} h_z(s) e^{n \phi (s)} ds.
$$
Before formulating a result on the large $n$ behavior of $F_n(z)$
we specify the type of contours $\gamma$ to be considered.
\begin{definition}
\label{def2.1}
A contour $\gamma$ is said to be \emph{admissible} if
\begin{itemize}
\item[(i)]
$\gamma$ is a smooth Jordan curve winding counterclockwise around the
origin, and
\item[(ii)]
$\gamma$ has a positive distance from the solid curve in Figure \ref{f2.1},
except for a part that lies in some neighborhood $U$ of $z_0=1$. In this
set $U$ the contour $\gamma$ coincides with the path of steepest descent
(dotted line in Figure \ref{f2.1}).
\end{itemize}
\end{definition}
The steepest descent analysis described above then gives
\begin{theorem}
\label{th2.1}
For any admissible contour $\gamma$ and $z \in \mathbb C \setminus \gamma$
we have
$$
F_n(z) = \frac{1}{\sqrt{2 \pi n}(1-z)} \left( 1 + O (n^{-1}) \right) \, ,
$$
where the error term is uniform for $z \in \mathbb C \setminus
(U_{\epsilon}(1) \cup \gamma)$ for any $\epsilon > 0$.
\end{theorem}
\begin{proof}
It remains to discuss the uniformity of the error bound. The only
step in the derivation of the asymptotic formula for which the uniformity
appears problematic is the estimate used in (\ref{eq2.20})
for points $z$
that lie arbitrarily close to $\gamma$ because
of the $(s-z)^{-1}$ term in the integrand.
However, since $z$ stays away some fixed distance $\epsilon$ from the
critical point $z_0$, one may deform the contour of integration for each $z$
in such a way that the value of $F_n$ does not change and such that
the deformed contour is separated from both, the solid line in Figure \ref{f2.1}
and from $z$ by a minimal distance that only depends on $\epsilon$.
\end{proof}

A somewhat more refined analysis provides uniform error bounds also
outside shrinking disks $|z-1| > n^{-\alpha}$ for any fixed $0 < \alpha
< 1/2$. In this case the error term $O(n^{-1})$ in Theorem \ref{th2.1}
has to be replaced by $O(n^{-1}|1-z|^{-2})$. For values of $z$ that lie
in shrinking discs $|z-1| = O(n^{-\alpha})$ it is still possible but
cumbersome to use the method described above for computing the asymptotics of $F_n(z)$. It is
precisely this situation in which we want to demonstrate how
ideas from the asymptotic analysis of Riemann-Hilbert problems can be used.

\section{Riemann-Hilbert analysis}
\label{S3}

In this section we utilize that the representation (\ref{eq1.20}), (\ref{eq1.40}) of $p_{n-1}$
is of Cauchy-type which enables us to employ Riemann-Hilbert techniques. The key for this is
that the function $F_n$ defined in (\ref{eq1.20}) is the Cauchy transform of
$e^{n\phi}$ with respect to the contour $\gamma$. Hence, and this is the
feature of the Cauchy transform that provides the link to Riemann-Hilbert problems,
$F_n$ is analytic in $\mathbb C \setminus \gamma$ and the values of
$F_n(z)$ differ by $e^{n\phi(s)}$ as $z$ approaches $s \in \gamma$ from opposite sides.
Moreover, these properties of $F_n$ together with the behavior of $F_n(z)$ as $z \to \infty$
characterize $F_n$ uniquely. More precisely, we will show in Lemma \ref{l3.1} below that
$F_n$ is the unique solution of the following scalar Riemann-Hilbert problem
($RHP$)$_1$:

\emph{Given an admissible contour $\gamma$ and $n \in \mathbb N$. Seek an
analytic function $Y \colon \mathbb C \setminus \gamma \to \mathbb C$ such that
\begin{itemize}
\item[(i)]
$Y_+(s) = Y_-(s) + e^{n\phi(s)}$ for $s \in \gamma$ \, ,
\item[(ii)]
$Y(z) \to 0$ for $|z| \to \infty$ \, .
\end{itemize}
}
Condition (i) is shorthand notation for the requirement that $Y$ has
continuous extensions from the interior of $\gamma$
(respectively from the exterior of $\gamma$) onto $\gamma$
which are denoted by $Y_+$ (respectively $Y_-$) and which satisfy relation (i).
As mentioned above the question of existence and uniqueness of solutions
for this Riemann-Hilbert problem is answered by the following

\begin{figure}
\includegraphics{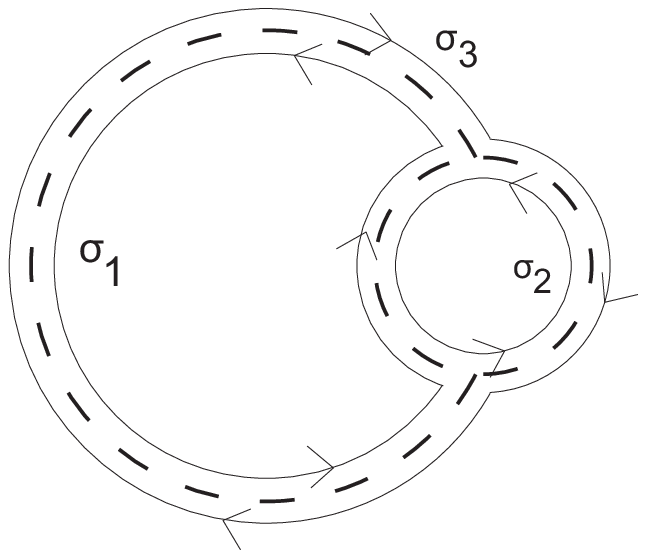}
\caption{}
\label{f3.1}
\end{figure}

\begin{lemma}
\label{l3.1}
$F_n$ as defined in (\ref{eq1.20}) is the unique solution of ($RHP$)$_1$.
\end{lemma}
\begin{proof}
Using (\ref{eq1.40}) it is easy to verify that $F_n$ indeed
solves ($RHP$)$_1$. In order to prove uniqueness one first shows that the
difference $\Delta$ of two solutions of ($RHP$)$_1$ is continuous across
$\gamma$ and hence entire. Liouville's theorem together with condition (ii)
then implies $\Delta(z) = 0$ for all $z \in \mathbb C$.
\end{proof}
Next we apply the method of constructing a local parametrix for the
Riemann-Hilbert problem in order to resolve the difficulties at the critical point
$z_0$.
We recall from Section \ref{S2} that the change of variables $z \to t$ with
$z= \lambda(it/\sqrt{n})$, $z \in U_0$, maps the contour
of steepest descent into the real axis and transforms $n \phi$ to normal form,
$n \phi(z) = - t^2$. This motivates the definitions
\begin{eqnarray}
\label{eq3.10}
h(\zeta) &:=& \frac{1}{2 \pi i} \int_{\mathbb R} \frac{e^{-u^2}}{u - \zeta} du\, ,
\quad \zeta \in \mathbb C \setminus \mathbb R \, , \\
P_n(z) &:=& h(-i\sqrt{n} \lambda^{-1}(z)) \, , \quad z \in U_0 \setminus \gamma\, .
\label{eq3.30}
\end{eqnarray}
We have $h_+(t)=h_-(t) + e^{-t^2}$ for $t\in \mathbb R$, where $h_{\pm}(t)$
denotes $\lim_{\eta \to 0+} h(t \pm i \eta$). This relation is a consequence of
standard properties of the Cauchy transform. Alternatively, since $h$ is the Cauchy
transform of an analytic function $e^{-u^2}$, one may also derive this relation
using only the calculus of residues. Substituting  $s = \lambda(it/\sqrt{n})$
we obtain
$$
(P_n)_+(s) = (P_n)_-(s) + e^{n \phi (s)} \quad \mbox{ for all } s \in U_0 \cap \gamma \, .
$$
$P_n$ is thus a local solution of the Riemann-Hilbert problem ($RHP$)$_1$ in $U_0$.
Moreover, $P_n$ is of a rather explicit nature, since $h$ is related to the well
studied complementary error function \cite{FCC} (see also (\ref{eq4.90}) below) and
the Taylor coefficients of $\lambda$ can be computed explicitly at $\xi=0$ to all
orders from the defining relation (\ref{eq2.30}). Note, however, that this local
solution $P_n$ cannot be continued to the global solution because $\lambda^{-1}(z)$ has
singularities outside of $U_0$.
The procedure to take advantage of this local parametrix is strikingly simple. Choose
$\epsilon > 0$ such that the closed disc $\overline{U_{2\epsilon}(1)}$ is contained in $U_0$.
We then set
\begin{equation}
\label{eq3.60}
\tilde{m}(z) := \left\{
\begin{array}{ll} Y(z),&\mbox{ for } z \in \mathbb C \setminus ( \gamma \cup \overline{U_{2 \epsilon} (1)})\\
Y(z) - P_n(z),&\mbox{ for } z \in U_{2 \epsilon} (1) \setminus \gamma
\end{array}
\right.
\end{equation}
where $Y=F_n$ denotes the unique solution of ($RHP$)$_1$. Observe that
$\tilde{m}_+(s) = \tilde{m}_-(s)$  for $s \in \gamma \cap U_{2 \epsilon}(1)$
since the jumps of $Y$
and $P_n$ cancel each other. We may therefore extend $\tilde{m}$
within $U_{2 \epsilon}(1)$ to an analytic function. We denote this function with a slightly extended
domain of definition by $m$. It is obvious that $m$ again solves a Riemann-Hilbert problem. More precisely,
denote
$$
\Gamma_1 := \partial U_{2 \epsilon}(1) \, , \quad
\Gamma_2 := \gamma \setminus U_{2 \epsilon}(1)\, , \quad
\Gamma := \Gamma_1 \cup \Gamma_2 \, .
$$
The dashed line in Figure \ref{f3.1} provides a sketch of $\Gamma$.
It is straightforward to verify that the function $m$ is a solution of ($RHP$)$_2$:

\emph{Seek an analytic function $M:\mathbb C \setminus \Gamma \to \mathbb C$ such that
\begin{itemize}
\item[(i)]
$M_+(s) = M_-(s) - P_n(s)$ for $s \in \Gamma_1 \setminus \Gamma_2$\\
$M_+(s) = M_-(s) + e^{n\phi(s)}$ for $s \in \Gamma_2$
\item[(ii)]
$M(z) \to 0$ for $|z| \to \infty$ \, .
\end{itemize}
}
As above $M_{\pm}$ denote the continuous extensions of $M$ from the
interior ($+$) and exterior ($-$).
We show in Remark \ref{r3.1} below that again we can use the Cauchy transform
to solve the Riemann-Hilbert problem ($RHP$)$_2$:
\begin{equation}
\label{eq3.50}
m(z) = \frac{1}{2 \pi i} \int_{\Gamma} \frac{m_+(s)-m_-(s)}{s-z}ds =
\frac{1}{2 \pi i} \int_{\Gamma_1} \frac{-P_n(s)}{s-z}ds +
\frac{1}{2 \pi i} \int_{\Gamma_2} \frac{e^{n \phi (s)}}{s-z}ds\, .
\end{equation}
The orientations of $\Gamma_1$ and $\Gamma_2$ are chosen to be counterclockwise
so that the $+$ side always lies to the left of the contour.
Observe that all $z \in U_{\epsilon}(1)$ have  at least distance $\epsilon$
from $\Gamma$ so that the contour of integration has been moved away from
the singularity $(s-z)^{-1}$
of the integrand.
Moreover, the integral over $\Gamma_2$ is
exponentially small so that the asymptotic expansion of $m$ is solely determined
by the integral over $\Gamma_1$. Before we formulate our theorem on the asymptotics
of $F_n(z)$ for $z$ near $z_0=1$ we present an elementary proof of (\ref{eq3.50})
which does not make use of the fact that $m$ is the unique solution of ($RHP$)$_2$.
Nevertheless, ($RHP$)$_2$ provides a clear explanation why (\ref{eq3.50}) holds:
The solution of this scalar Riemann-Hilbert problem is simply given by the Cauchy
transform of the jump $m_+ - m_-$on $\Gamma$.

\begin{remark}
\label{r3.1}
One may show the first equality of (\ref{eq3.50}) by applying the calculus of
residues to
$$
\int_{\sigma_1} \frac{m(s)}{s-z} ds +
\int_{\sigma_2} \frac{m(s)}{s-z} ds +
\int_{\sigma_3} \frac{m(s)}{s-z} ds
$$
where $\sigma_1$, $\sigma_2$ and $\sigma_3$ denote the closed curves shown in
Figure \ref{f3.1} and by taking the limit as these curves approach $\Gamma$
(dashed line). Note that condition (ii) in ($RHP$)$_2$ implies that $m(s)/(s-z)$
has a vanishing residue at $s= \infty$.
\end{remark}
We are now ready to state our main result.
\begin{theorem}
\label{th3.1}
There exists $\epsilon > 0$ and functions $g_j$ analytic in $U_{\epsilon}(1)$
such that for any admissible curve $\gamma$ and any $r \in \mathbb N$ we have
\begin{equation}
\label{eq3.70}
F_n(z) = P_n(z) + \frac{1}{\sqrt{2 \pi n}}
\left(\sum_{j=0}^{r-1} \frac{g_j(z)}{n^j} + O \left(\frac{1}{n^r} \right)
\right) \, ,
\end{equation}
where the error term is uniform for $z\in U_{\epsilon}(1)\setminus\gamma$
and the Taylor coefficients of all $g_j$ are explicitly computable at $z_0=1$,
e.g.
$$
g_0(z) = \frac{1}{3} - \frac{1}{12}(z-1) + O\left((z-1)^2\right)\, .
$$
\end{theorem}
\begin{proof}
Choose $\epsilon > 0$ such that the closed disc $\overline{U_{2\epsilon}(1)}$
is contained in $U_0 \cap U$ where $U_0$ is defined above (\ref{eq2.30}) and $U$ is
determined by the curve $\gamma$ (cf. Definition \ref{def2.1}(ii)).
From the discussion of the present section (see in particular (\ref{eq3.60}), (\ref{eq3.50}))
it follows that the unique solution $F_n = Y$ of
($RHP$)$_1$ can be written for $z \in U_{\epsilon}(1)$ in the form
$$
F_n(z) = P_n(z) + m(z) = P_n(z) - \frac{1}{2 \pi i} \int_{\Gamma_1} \frac{P_n(s)}{s-z}ds
+ O(e^{-nc})\, .
$$
Observe that we may not use the calculus of residues to evaluate
$\int_{\Gamma_1} \frac{P_n(s)}{s-z}ds$ since $P_n$ is not a meromorphic
function. Nevertheless, using
$$
\frac{1}{u-\zeta} =
-\sum_{j=0}^{2r} \frac{u^j}{\zeta^{j+1}} + \frac{u^{2r+1}}{\zeta^{2r+1}(u-\zeta)}
$$
together with the fact that the Cauchy transform
$\frac{1}{2 \pi i} \int_{\mathbb R} \frac{s^{2r+1}e^{-s^2}}{s-z}ds$ is
bounded on $\mathbb C \setminus \mathbb R$ (which can e.g. be proved using
a contour deformation argument)
one derives that for each $r \in \mathbb N$ there exists a constant $C$ such that
$$
\left|
h(\zeta) + \frac{1}{2 \pi i \zeta} \sum_{j=1}^{r} \frac{\Gamma(j+\frac{1}{2})}{\zeta^{2j}}
\right| \leq
\frac{C}{|\zeta|^{2r+1}} \,
$$
for all $\zeta \in \mathbb C \setminus \mathbb R$.
From the definition of $P_n$ in terms of $h$ we obtain
\begin{equation}
\label{eq3.80}
P_n(s) = \frac{1}{\sqrt{2 \pi n}} \left(
\sum_{j=0}^{r} \frac{G_j(s)}{n^j(s-1)^{2j+1}} + \mbox{rest }
\right)\, .
\end{equation}
By choosing $\epsilon$ smaller, if necessary, we can ensure that all $G_j$
are analytic and non-zero in some open neighborhood of $\overline{U_{2\epsilon}(1)}$.
The term ``rest'' in (\ref{eq3.80}) is
of order $n^{-r}$ uniformly for $s \in \partial U_{2 \epsilon} (1)$.
The calculus of residues then leads to the desired representation of $F_n$
and the Taylor coefficients of $g_j$ can be derived from the Taylor coefficients
of $\lambda^{-1}$ by explicit computation.
\end{proof}
\begin{remark}
\label{r3.2}
Replacing $U_{2 \epsilon}(1)$ by $U_{\epsilon /2}(1)$ in (\ref{eq3.60})
and adapting the arguments above one obtains
an expansion of $F_n(z)$ also for $|z-1| > \epsilon$.
Indeed, one can prove that
there exist polynomials $h_j$ of degree $2j$ such that for all $r \in \mathbb N$
we have
\begin{equation}
\label{eq3.90}
F_n(z) = \frac{1}{\sqrt{2 \pi n}(1-z)} \left( 1+
\sum_{j=1}^{r-1} \frac{h_j(z)}{n^j (z-1)^{2j}} + O\left(\frac{1}{n^{r}} \right)
\right) \, ,
\end{equation}
where the error term is uniform for $\epsilon \leq |z-1|\leq 2$. Again the polynomials
$h_j$ may be computed explicitly, e.g. $h_1(z) = -\frac{1}{12}(z^2+10z+1)$. This result
sharpens the statement of Theorem \ref{th2.1}. It is also not difficult to see that
the expansion (\ref{eq3.90}) holds outside shrinking circles $n^{-\alpha} \leq |z-1|
\leq 2$ for any fixed $0 < \alpha < 1/2$. In this case the error term has to be
replaced by $O\left(\frac{1}{n^{r}|1-z|^{2r}} \right)$ to ensure uniformity in $z$.
\end{remark}

\section{Locating the zeros of the exponential sum}
\label{S4}

The arguments and results put forward in this section will be presented
in detail in a later publication, where we will construct for a
specific admissible contour $\gamma$ precisely $n-1$ different solutions
(together with their asymptotic expansions) of the equation
\begin{equation}
\label{eq4.10}
e^{n \phi (z)} = F_n (z)
\end{equation}
in the interior of
$\gamma$, provided $n$ is sufficiently large.
Since $p_{n-1}(nz)$ has at most $n-1$ different zeros (\ref{eq1.40}) implies that
we have located all the zeros of $p_{n-1}$. We will find the
solutions of (\ref{eq4.10}) by showing for each $1 \le k \le n-1$
that the equation
\begin{equation}
\label{eq4.20}
G_n(z) := \tilde{\phi}(z) - \frac{1}{n} \ln F_n(z) = - \frac{2 \pi i k}{n}
\end{equation}
has one solution in the unit disc $U_1(0)$.
Here $\tilde{\phi}$ is defined as $\phi$ in (\ref{eq1.20})
with the only difference that the branch of the logarithm
is now chosen in such
a way that its imaginary part takes values in $(0, 2 \pi)$ rather than in
$(- \pi, \pi)$. Existence and asymptotic expansions of the roots of
(\ref{eq4.20}) are obtained by a standard procedure. First one constructs
solutions $\alpha_{k,n}$ of
\begin{equation}
\label{eq4.30}
A_n(z) = - \frac{2 \pi i k}{n}
\end{equation}
for some approximation $A_n$ of $G_n$. Then a contraction mapping argument
will be used to conclude that the original equation (\ref{eq4.20})
has a solution $z_{k,n}$ close to $\alpha_{k,n}$. Since we can prove uniform bounds on
the derivatives of the error terms in (\ref{eq3.70}), (\ref{eq3.90}) we obtain in addition
an asymptotic expansion for $z_{k,n}$ in terms of $\alpha_{k, n}$.

\subsection{Zeros away from the critical point}
\label{S4.1}
In this case it suffices to use the crude approximation
$$
A_n(z) = \tilde{\phi}(z) = z - 1 - \widetilde{\ln}(z)\, ,
$$
where $\widetilde{\ln}$ denotes the branch of the logarithm described above. The solutions
$\alpha_{k,n}$ of (\ref{eq4.30}) all lie on the Szeg\H{o} curve $D_{\infty}$.
One can show that the distance of $\alpha_{k,n}$ and $z_{k, n}$ is of the same order as
$A_n - G_n = \frac{1}{n}(\ln F_n)=O((\ln n)/n)$ (see Theorem \ref{th2.1}).
Note that this result only holds for $z_{k,n}$ that lie in a compact subset of $\mathbb C
\setminus \{ 1 \}$. An extended version of Theorem \ref{th2.1} as discussed at the end of
Section \ref{S2} allows to
include also those solutions
$z_{k, n}$ that lie outside shrinking discs $|z-1| \geq n^{-\alpha}$ with $0 < \alpha < 1/2$. The
distance between $\alpha_{k,n}$ and $z_{k, n}$ is then of the order $\frac{\ln n}{n|1-z_{k,n}|^2}$.
We state our result on the asymptotic expansion of $z_{k, n}$.
\begin{theorem}
\label{th4.2}
There exist polynomials $Q_j(x, y)$ of degree $j$ in the variable $y$ and of degree
$\leq 2j-2$ in the variable $x$ such that for $0 < \beta < 1$,
$n^{\beta} < k \leq n/2$, and $r \in \mathbb N$ we have
\begin{equation*}
     z_{k,n}
     = \alpha_{k,n}\left(1+\sum^{r-1}_{j=1}
     \frac{Q_j\left(\alpha_{k,n},\ln \left[\sqrt{2\pi n}(\alpha_{k,n}-1)\right]\right)}
     {n^j(1-\alpha_{k,n})^{2j-1}}\right) +
     O \left[
      \left(\frac{\ln n}{n}\right)^r\left(\frac{n}{k}\right)^{r-\frac{1}{2}}\right] \, ,
\end{equation*}
where the constant in the error term only depends on the choice of $\beta$ and $r$.
The polynomials $Q_j$ can be computed explicitly. For example, we have
\begin{eqnarray*}
     Q_1(x,y)=-\frac{1}{2}y;\quad Q_2(x,y)=-\frac{1}{8}y^2+\frac{1}{2} xy-\frac{1}{12} (x^2+10x+1).
\end{eqnarray*}
\end{theorem}
For $r=1$ (i.e. without correction terms) this result was first proved in \cite[(A.47)]{BM}.

Faster rates of convergence
can be achieved by using a better approximation of $F_n$ in the definition of
$A_n$ that is provided by Theorem \ref{th2.1}, namely
$$
A_n(z) = \tilde{\phi}(z) + \frac{1}{n} \ln (\sqrt{2 \pi n} (1-z))\, .
$$
The faster rate of approximation comes at the price that the corresponding
approximate solutions $\alpha_{k,n}$ now lie on $n$ -- dependent curves $D_n$ rather
than on the Szeg\H{o} curve. We leave the corresponding expansion of the zeros in terms
of such approximate solutions $\alpha_{k, n}$ for a later publication. Note that
both \cite[(A.48)]{BM} and \cite{CVW} also work with such better approximations.

\subsection{Zeros near the critical point}
\label{S4.2}
In order to formulate the result we first need to recall the definition of the
complementary error function
\begin{eqnarray*}
     \mbox{erfc}\, (z)=1-\frac{2}{\sqrt{\pi}}\int\limits^z_0e^{-t^2}\, dt=
     \frac{2}{\sqrt{\pi}}\int\limits^{\infty}_ze^{-t^2}\, dt\quad \textrm{for}\  z\in \mathbb C,
\end{eqnarray*}
where the path of integration of the latter integral is subject to the restriction
arg($t$) $\to \alpha$ with $|\alpha| < \frac{\pi}{4}$ as $t \to \infty$ along the path.
It is well known (see e.g. \cite{FCC}) that all the zeros of this function lie in the
second and third quadrant of the complex plane
(i.e. in the regions $\frac{\pi}{2}<\arg (z) <\pi$ and $-\pi < \arg(z)<-\frac{\pi}{2}$)
and that
in the second quadrant there are countably many zeros of the complementary error function.
We denote these zeros
by $w_k,k\in \mathbb N$, and we can order them by modulus $|w_k|<|w_{k+1}|$.
Our result on the
solutions $z_{k, n}$ of (\ref{eq4.20}) reads as follows.
\begin{theorem}
\label{th4.1}
There exist polynomials $q_j$ of degree $j$ such that for $0 < \beta < 1$,
$1 < k < n^{\beta}$, and $r \in \mathbb N$ we have
\begin{equation*}
     z_{k,n} = 1+\sum^{r-1}_{j=1} \frac{q_j(\sqrt{2}w_k)}{n^{j/2}} +
     O\left(\left(\frac{k}{n}\right)^{r/2}\right) \, ,
\end{equation*}
where the constant in the error term only depends on the choice of $r$ and $\beta$. Moreover,
the polynomials $q_j$ may be computed explicitly, e.g.
\begin{eqnarray*}
     q_1(x)=x;\quad q_2(x)=\frac{x^2-1}{3};\quad q_3(x)=\frac{x^3-7x}{36}.
\end{eqnarray*}
\end{theorem}
Such a result was proved for $r=2$ and $\beta < 1/3$ in \cite[(A.34)]{BM} sharpening and extending previous
results of \cite{NR}, \cite{CVW}.
Since $|w_k| \sim \sqrt{2 \pi k}$ for $k \to \infty$ it follows that the $j$-th term in the expansion
above is of order $(\frac{k}{n})^{j/2}$. In particular, we have that $|z_{k, n} - 1|$
is of order $\sqrt{\frac{k}{n}}$.
Consequently, the expansion of $z_{k,n}$ in terms of the zeros of the complementary error function
holds for all zeros in shrinking circles of size $n^{(\beta - 1)/2}$. Thus, for any $\varepsilon > 0$ and
$n^{-(1/2)+\varepsilon} < |z_{k, n} -1| < n^{-\varepsilon}$ subsections \ref{S4.1} and \ref{S4.2} provide different
expansions for the
solutions of (\ref{eq4.20}). A short calculation shows that Theorem \ref{th4.1} yields better approximations
in shrinking circles of size
$O(n^{-1/3})$, otherwise the approximation with $\alpha_{k,n}$ on the Szeg\H{o} curve
(see Theorem \ref{th4.2}) is more advantageous.

We finish by explaining how the complementary error function enters the picture in the proof of
Theorem \ref{th4.1}.
Introducing the auxiliary function
$$
v(\zeta) := e^{\zeta^2} \mbox{ erfc}(\zeta),
$$
one verifies that $v(\zeta) = 2 e^{\zeta^2} - v(-\zeta)$ holds for $\zeta \in \mathbb R$ and hence
by the identity principle on all of the complex plane. Since $v$ and the complementary error function
have the same set of roots we obtain $2e^{w_k^2} v(-w_k)^{-1}=1$.
Elementary estimates on the arguments show that the correct value of the logarithm is given by
\begin{equation}
\label{eq4.80}
w_k^2 - \ln(v(-w_k)/2) + 2 \pi i k = 0 \, , \quad \mbox{ for all } k \in \mathbb N \, .
\end{equation}
Next we state the relation between $v$ and the function $h$ defined in (\ref{eq3.10}).
For all $\zeta$ with positive real part the following relation holds
\begin{equation}
\label{eq4.90}
v(\zeta) =2 h(i \zeta) \, .
\end{equation}
To see this one verifies that
$g(\zeta):= 2 h(i \zeta)e^{-\zeta^2}$ has the same derivative and the same limiting
behavior for $\zeta \to \infty$ as the complementary error function. Recall
that the function $h$ was used to define the parametrix $P_n$ (\ref{eq3.30}). Following the procedure
described at the beginning of the present section and keeping the result of Theorem \ref{th3.1}
in mind we approximate
$G_n$ by
$$
A_n(z) := \tilde{\phi}(z) - \frac{1}{n} \ln P_n(z) \, .
$$
Setting $\alpha_{k, n} := \lambda(w_k/\sqrt{n})$, using that $\phi$ and $\tilde{\phi}$
agree on the upper half plane, and that the real part of $-w_k$ is positive,
we obtain from (\ref{eq2.30}), (\ref{eq4.90}), and (\ref{eq4.80})
$$
A_n(\alpha_{k,n}) = \phi(\lambda(w_k/\sqrt{n})) - \frac{1}{n} \ln h(-i w_k) =
\frac{1}{n}[ w_k^2 - \ln (v(-w_k)/2) ] = - \frac{2 \pi i k}{n}
$$
satisfying (\ref{eq4.30}) as desired.

\bibliographystyle{amsplain}

\begin{thebibliography}{99}

\bibitem{ACV} V. V. Andrievskii, A. J. Carpenter, and R. S. Varga,
\textit{Angular distribution of zeros of the partial sums of $e\sp z$ via the
solution of inverse logarithmic potential problem},
Comput. Methods Funct. Theory {\bf 6} (2006), no. 2, 447--458.

\bibitem {BM} P. Bleher and R. Mallison,
\textit{Zeros of sections of exponential sums}, International Mathematics
Research Notices {\bf 2006}, Art. ID 38937, 49 pages.

\bibitem {B} J. D. Buckholtz, \textit{A characterization of the exponential series}, The American Mathematical Monthly {\bf 73} (1966), no. 4, part II, 121--123.

\bibitem {CVW} A. J. Carpenter, R. S. Varga, and J. Waldvogel,
\textit{Asymptotics for the zeros of the partial sums of $e^{z}$. $I$}, The Rocky Mountain Journal of Mathematics {\bf 21} (1991), no. 1, 99--120.

\bibitem {DZ1} P. Deift and X. Zhou,
\textit{A steepest descent method for oscillatory Riemann-Hilbert problems. Asymptotics
for the mKdV equation}, Ann. of Math. (2) {\bf 137} (1993), no. 2, 295--370.

\bibitem {DZ2} P. Deift and X. Zhou,
\textit{Asymptotics for the {P}ainlev\'{e} {II} equation}, Comm. Pure Appl.
Math. {\bf 48} (1995), no. 3, 277--337.

\bibitem {FCC} H. E. Fettis, J. C. Caslin, and K. R. Cramer,
\textit{Complex zeros of the error function and of the complementary error function},
Mathematics of Computation {\bf 27} (1973), no. 122, 401--407.

\bibitem{Miller} P. D. Miller,
\textit{Applied asymptotic analysis}, Graduate Studies in Mathematics {\bf 75},
Amer. Math. Soc., Providence, RI, 2006.

\bibitem{NR} D. J. Newman and T. J. Rivlin,
\textit{The zeros of the partial sums of the exponential function},
Journal of Approximation Theory {\bf 5} (1972), no. 4, 405--412,
\textit{Correction}: Journal of Approximation Theory {\bf 16} (1976), 299--300.

\bibitem {O} I. V. Ostrovskii,
\textit{On zero distribution of sections and tails of power series},
Entire Functions in Modern Analysis (Tel-Aviv, 1997), Israel Math. Conf. Proc., vol. {\bf 15},
Bar-Ilan University, Ramat Gan, 2001, pp. 297--310.



\bibitem{S} G. Szeg\H{o},
\textit{\"Uber eine Eigenschaft der Exponentialreihe},
Sitzungsberichte, Berliner Mathematische Gesellschaft {\bf 23} (1924), 50--64.

\bibitem{V} R. S. Varga,
\textit{Scientific computation on mathematical problems and conjectures},
CBMS-NSF Regional Conference Series in Applied Mathematics {\bf 60},
SIAM, Philadelphia, PA, 1990.

\bibitem{VC1} R. S. Varga and A. J. Carpenter,
\textit{Zeros of the partial sums of $\cos (z)$  and $\sin (z)$. $I$},
Numerical Algorithms {\bf 25} (2000), no. 1--4, 363--375.


\bibitem{Z} S. M. Zemyan,
\textit{On the zeroes of the $N$- th partial sum of the exponential series},
The American Mathematical Monthly {\bf 112} (2005), no. 10, 891--909.

















\end{thebibliography}

\end{document}